\documentclass[UTF8,9pt,a4paper]{amsart} 
\usepackage[utf8]{inputenc}

\usepackage[british]{babel}

\usepackage{footmisc}

 \usepackage[%
     style=phys,%
     biblabel=brackets,%
 ]{biblatex}
\renewbibmacro{in:}{}
\usepackage{url}
\usepackage[colorlinks]{hyperref}

\usepackage{graphicx}

\usepackage{csquotes}

\usepackage{amsthm}
\usepackage{amsmath,amssymb}
\usepackage{xfrac}
\usepackage{braket}

\usepackage{hyperref}
\usepackage{cleveref} 

\theoremstyle{remark}
\newtheorem{remark}{Remark}

\usepackage{epsfig}  		
\usepackage{epic,eepic}

\makeatletter
\def\paragraph{\@startsection{paragraph}{4}%
  \z@\z@{-\fontdimen2\font}%
  {\normalfont\bfseries}}
\makeatother

\newenvironment{nouppercase}{%
  \renewcommand{\uppercasenonmath}[1]{}}{}

\begin{document}

\begin{nouppercase}
\title[Stunned by Sleeping Beauty]{\large Stunned by Sleeping Beauty: \\
\normalsize How Prince Probability updates his forecast upon their fateful encounter}
\end{nouppercase}

\author{Laurens Walleghem}


\date{ \today \newline 
\href{mailto:laurens.walleghem@york.ac.uk}{laurens.walleghem@york.ac.uk} \newline Department of Mathematics, University of York, York, United Kingdom \newline INL -- International Iberian Nanotechnology Laboratory, Braga, Portugal }

\begin{abstract} 
    The Sleeping Beauty problem is a puzzle in probability theory that has gained much attention since Elga's discussion of it [Elga, Adam, Analysis \textbf{60} (2), p.143-147 (2000)].
    Sleeping Beauty is put asleep, and a coin is tossed. 
    If the outcome of the coin toss is Tails, Sleeping Beauty is woken up on Monday, put asleep again and woken up again on Tuesday (with no recollection of having woken up on Monday). 
    If the outcome is Heads, Sleeping Beauty is woken up on Monday only. 
    Each time Sleeping Beauty is woken up, she is asked what her belief is that the outcome was Heads.     
    What should Sleeping Beauty reply?
    In literature arguments have been given for both 1/3 and 1/2 as the correct answer.
    In this short note we argue using simple Bayesian probability theory why 1/3 is the right answer, and not 1/2. 
    Briefly, when Sleeping Beauty awakens, her being awake is nontrivial extra information that leads her to update her beliefs about Heads to 1/3. 
    We strengthen our claim by considering an additional observer, Prince Probability, who may or may not meet Sleeping Beauty. 
    If he meets Sleeping Beauty while she is awake, he lowers his credence in Heads to 1/3.
    We also briefly consider the credence in Heads of a Sleeping Beauty who knows that she is dreaming (and thus asleep).
\end{abstract}

\begingroup
\def\uppercasenonmath#1{} 
\maketitle
\endgroup

\tableofcontents


\newpage

\section{Introduction}
\paragraph{The Sleeping Beauty problem} We follow \cite[Appendix C]{jones2024thinking} for the description of the Sleeping Beauty\footnote{We will sometimes abbreviate Sleeping Beauty as SB or simply Beauty.} problem\footnote{The SB problem was first written down by Elga in \cite{elga2020self}, who attributed it to Robert Stalnaker and Arnold Zuboff. Similar problems were posed by Aumann, Hart and Perry in \cite{aumann1997forgetful} and Piccione and Rubinstein in \cite{piccione1997interpretation}.}.
\begin{quote}
    \textit{The Sleeping Beauty problem: }
    Sleeping Beauty participates in an experiment in which she is put to sleep on Sunday, and a fair coin is tossed. If the outcome is Heads, she will be woken up on Monday only. If the outcome is Tails, she will be woken up on Monday and Tuesday. More specifically, for outcome Tails Sleeping Beauty is woken up on Monday and put back to sleep with her memory of her awakening erased, to be woken up on Tuesday again. 
    Sleeping Beauty cannot distinguish between any of these awakenings. 
    When woken up, Sleeping Beauty is asked for her credence that the outcome of the coin toss was Heads. 
\end{quote}

As a proposal to make the question a bit more clear, we may ask Sleeping Beauty, when she is woken up, the following question: `What do you now think the outcome was, Heads or Tails? And what is your credence now for Heads?'

\vspace{0.5cm}

Two possibilities for Sleeping Beauty to answer have been debated heavily in literature. The Halfers argue that 1/2 is the correct answer, the main argument being the absence of any new information for Beauty when she wakes up so that she should stick to her credence 1/2 she has on Sunday, whereas the Thirders argue that 1/3 is the correct answer.
Many works in favour of Elga's argument for 1/3 \cite{elga2020self}, and against Lewis's 1/2 arguments of \cite{lewis2001sleeping} have appeared \cite{elga2020self,elga2004defeating,hitchcock2004beauty,draper2008diachronic,arntzenius2003some,dorr2002sleeping,van2013sleeping,horgan2004sleeping,weintraub2004sleeping,monton2002sleeping,titelbaum2012embarrassment,titelbaum2013quitting,neal2006puzzles,stalnaker2010our}.
The arguments for 1/2 mostly claim that no extra information is present for Beauty when she is woken up.\footnote{The fact that this is false was already noted (among others) in \cite{weintraub2004sleeping,horgan2004sleeping,hitchcock2004beauty}.} 
Here we argue that this claim is wrong using simple Bayesian probability theory and a slight extension of the protocol, with an outsider Prince Probability who may or may not meet Sleeping Beauty. We also consider the credence of a Sleeping Beauty who knows that she is dreaming. We believe that the strength of our argument is its simplicity. 

\vspace{0.5cm}

\paragraph{Outline} 
We explain our main argument using simple Bayesian probability in \Cref{sec:SB_Thirder}, where we also describe how Prince Probability updates his credence  upon meeting Sleeping Beauty.
We consider some generalisations and variations of the SB problem in \Cref{sec:SB_variations}.
We relate to Elga's argument and literature in \Cref{sec:literature}.
Finally, in \Cref{sec:SB_and_FR} we comment on the relation of the SB problem to the Frauchiger--Renner paradox \cite{frauchiger2018quantum}, following \cite{jones2024thinking}. 

\vspace{0.5cm}

\section{A Thirder's Tale of Beauty and the Prince} \label{sec:SB_Thirder}

In this section we will argue how Sleeping Beauty knowing that she is awake is nontrivial useful information for Sleeping Beauty, leading to the Thirder position. First, we consider Beauty's point of view, but perhaps most convincing is the story including an outsider Prince Probability, who may or may not meet Beauty and updates his credence in Heads accordingly, which we discuss subsequently. 
Finally, we argue what credence Beauty would assign if she knew she were dreaming (and thus asleep), the same credence that Prince would assign upon meeting a sleeping Sleeping Beauty.

\vspace{0.5cm}

\paragraph{Beauty knowing that she is awake is nontrivial information.} In the Sleeping Beauty protocol, Beauty is asked to give her credence for Heads when being woken up. 
Importantly, here the fact that Beauty is awake is nontrivial: she could have been in a state of sleep too.
Suppose for example that each time Beauty is woken up she is awake for 1 hour. The event of her being awake under Heads then has probability \begin{equation}
    p(\text{SB awake}|H) = p(\text{SB awake $\land$ Monday}|H) = 1/(2\cdot 24) = 1/48,
\end{equation} as she is awake for 1 hour only in the two days (48 hours) of the experiment.\footnote{Here and further on, the events $H$ and $T$ stand for the outcome of the coin toss being Heads and Tails, respectively.}
On the other hand, under Tails we have \begin{equation}
\begin{split}
            p(\text{SB awake}|T) &= p(\text{SB awake $\land$ Monday}|T) + p(\text{SB awake $\land$ Tuesday}|T) \\ &= 1/48+1/48 = 1/24.
\end{split}
\end{equation} 
The probability of Beauty being awake is thus \begin{equation}
    p(\text{SB awake}) = p(H) p(\text{SB awake}|H) + p(T) p(\text{SB awake}|T) = \frac{1}{2} \cdot \frac{1}{48}+\frac{1}{2} \cdot \frac{1}{24} = \frac{1}{2} \cdot \frac{3}{48}.
\end{equation} 
Therefore, the credence Beauty gives to Heads while being awake is \begin{equation} \label{sec:Heads_prob_main}
    p(H|\text{SB awake}) = \frac{p(H) p(\text{SB awake}|H)}{p(\text{SB awake})} = \frac{1/2 \cdot 1/48}{1/2 \cdot 3/48} = 1/3.
\end{equation} 
Analogously, we find that $p(T|\text{SB awake})=2/3$. 
Thus indeed, while being awake Sleeping Beauty gives credence 1/3 to Heads, using the fact that she is awake, because she is more likely to be awake when the outcome is Heads. 

\vspace{0.5cm}

\begin{remark}[Generalisation to $N$ days and arbitrary `awake' time of Sleeping Beauty.]
    We can generalise this argument to any time span Beauty is awake for when woken up before being put to sleep again, and the experiment taking any number $N$ of days. Namely, consider the variation where for Tails Beauty is woken up for $N$ consecutive days, whereas for Heads only once in those $N$ days. Furthermore, each time she is woken up, Beauty is awake for $z$ hours that day. Then, for outcome Tails as Beauty is woken up $N$ consecutive days we have $P(\text{SB awake}|T) = (N \cdot z)/(24 \cdot N) = x$, where $x=z/24$.
    For Heads she is woken up only once during those $N$ days, so $P(\text{SB awake}|H) = (1 \cdot z)/(N \cdot 24) = \frac{x}{N}$.
    Therefore, 
    \begin{equation}
    \begin{split}
        P(\text{SB awake})&= p(H)p(\text{SB awake}|H) + p(T) p(\text{SB awake}|T)\\ & =  \frac{1}{2} \cdot \frac{x}{N} + \frac{1}{2} \cdot x = \frac{1}{2} \cdot x \cdot \frac{N+1}{N},        
    \end{split}
    \end{equation} and so 
    \begin{equation} \label{eq:many_awakenings_H}
        P(H|\text{SB awake}) = \frac{P(\text{SB awake}|H)p(H)}{p(\text{SB awake})} = \frac{1/2 \cdot x \cdot \frac{1}{N}}{1/2 \cdot x \cdot \frac{N+1}{N}} = \frac{1}{N+1},
    \end{equation} and 
    \begin{equation}
        P(T|\text{SB awake}) = \frac{P(\text{SB awake}|T)p(T)}{p(\text{SB awake})} = \frac{1/2 \cdot x}{1/2 \cdot x \cdot \frac{N+1}{N}} = \frac{N}{N+1}.
    \end{equation}
    Thus, if the experiment takes $N$ days, when woken up and asked for her belief in Heads, Beauty replies having credence $1/(N+1)$ for Heads.
    For example, if $N=9$, i.e. for Tails Beauty is woken up 9 consecutive days, each time she is woken up SB gives credence $1/10$ to Heads.
    As expected, this result is thus indeed independent of $x$, i.e., independent of the time Beauty is awake for each time she is woken up. (There is of course the restriction $x < 1$ that Beauty cannot be awake for a full day or more each time she is woken up.)
    \end{remark}

\vspace{0.5cm}

\begin{figure}[t]
         \centering
\includegraphics[width=0.8\columnwidth]{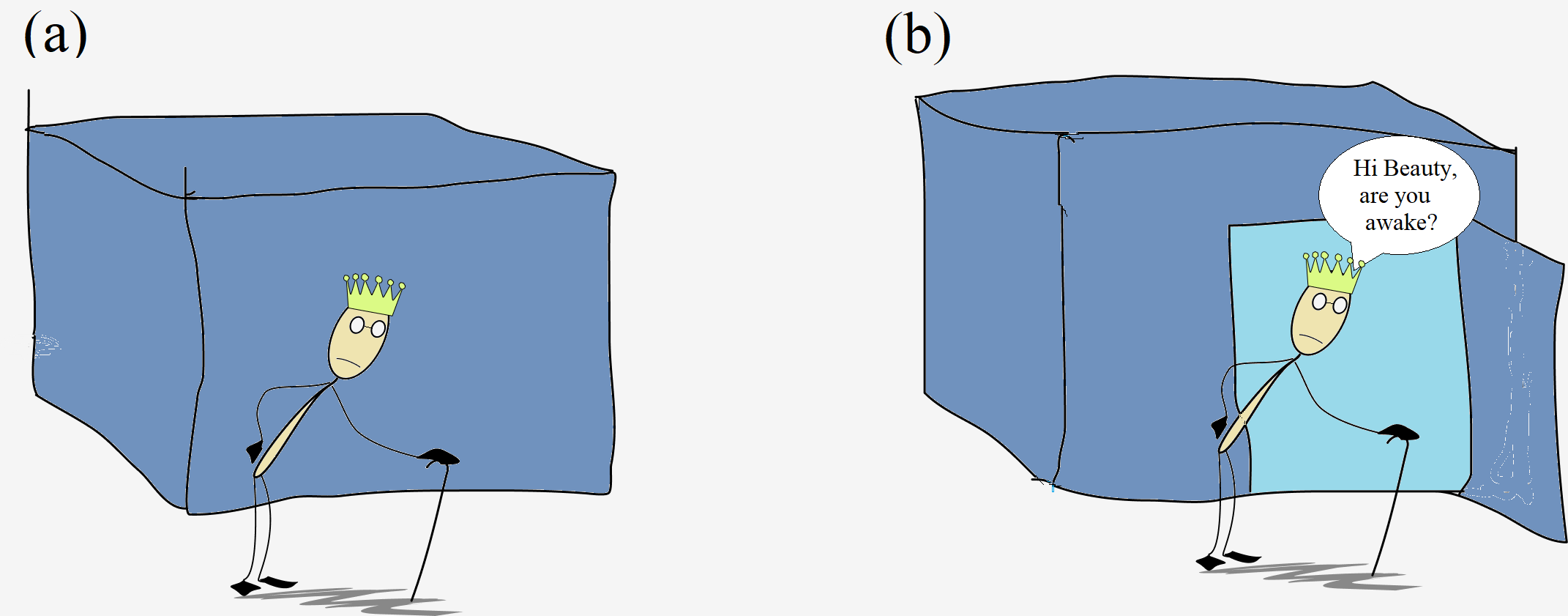}
         \caption{A schematic overview of the extension of the Sleeping Beauty problem that includes an outsider's perspective (Prince Probability's) as described in \Cref{sec:SB_Thirder}. (a): Sleeping Beauty resides in a sealed lab and an outsider Prince Probability assigns a probability 1/2 to Heads. (b): Prince Probability opens Sleeping Beauty's lab and can now see whether she is awake or asleep. If she is awake, Prince Probability lowers his credence for Heads to 1/3. If she is asleep, he increases his credence for Heads.
          }
         \label{fig:PrinceProbability_SB}
\end{figure}

\paragraph{Perspective of outsider Prince Probability.} 
Consider an outsider, Prince Probability\footnote{We will also refer to Prince Probability simply as Prince or PP.}, and for concreteness let Sleeping Beauty always reside in a sealed lab. When not partaking in the protocol, Prince Probability assigns credence 1/2 to Tails. Prince Probability has no idea of which day it is, and Prince Probability visits Beauty's lab. He sees Beauty either asleep or awake. If Prince Probability sees Beauty being awake, he will lower his credence to 1/3 for Heads because under Tails she is awake for more time. The calculation is analogous to the one right above in \cref{sec:Heads_prob_main}. 
Upon meeting the Prince, Beauty does not seem to get any relevant information for Heads versus Tails, suggesting that she should indeed always assign credence 1/3 to Heads when she is awake. We have sketched this protocol in \Cref{fig:PrinceProbability_SB}. 

\vspace{0.5cm}

\paragraph{What if Sleeping Beauty knew she was dreaming?} 
Consider the Sleeping Beauty problem, but now Beauty knows she is dreaming (and thus asleep, assuming she always dreams while asleep). Using similar Bayesian arguments as above, she will now increase her credence in Heads. 
Suppose Sleeping Beauty is awake for $x$ hours each time she is woken up (so $x < 24$ and we assume $x$ to be known to her). We have that $P(\text{SB asleep}|H) = (48-x)/48 = 1-y$ and $P(\text{SB asleep}|T) = (48-2x)/48=1-2y$ where we have defined $y=x/48$ (with now $y < 1/2$). 
Therefore, \begin{equation}
\begin{split}
    P(\text{SB asleep}) &= P(H)  P(\text{SB asleep}|H) + P(T)  P(\text{SB asleep}|T) \\ &= 1/2 \cdot (1-y) + 1/2 \cdot (1-2y) = 1/2 \cdot (2-3y),
\end{split}
\end{equation} and thus \begin{equation}
    P(H|\text{SB asleep}) = \frac{P(H) P(\text{SB asleep}|H)}{P(\text{SB asleep})} = \frac{1/2 \cdot (1-y)}{1/2 \cdot (2-3y)} = \frac{1-y}{2-3y}.
\end{equation} 
Thus, knowing that she is currently dreaming, Sleeping Beauty assigns credence $(1-y)/(2-3y)$ to Heads in this scenario.
This is also the credence an outsider Prince Probability would give to Heads when he sees that SB is asleep.
Importantly, through $y$ this result does depend on the time span Sleeping Beauty is awake for each time when woken up.
For example, if each time she wakes up, SB is awake for 23 (1) hours, she gives credence 0.926... (0.505...) to Heads while knowing that she is dreaming.
This difference and dependence on $y$ make sense intuitively: for $x=23$ Beauty is asleep for 25 hours if Heads and only 2 hours if Tails, whereas for $x=1$ Beauty is asleep for 47 hours if Heads and for 46 hours if Tails, a much smaller difference.  

\vspace{0.5cm}

\section{Variations of the Sleeping Beauty problem}
\label{sec:SB_variations}

In this section we discuss two variations of the SB problem, namely the copy-SB problem of \cite{jones2024thinking} and the case of an unfair coin toss. We have already discussed the case of many awakenings above in \cref{eq:many_awakenings_H}.

\subsection{The copy Sleeping Beauty version of \cite{jones2024thinking}}
\label{sec:copy_SB}

We now discuss a slightly altered version of the copy-SB version of \cite{jones2024thinking}, which has a resolution similar to the one we proposed above in \Cref{sec:SB_Thirder} for the original SB problem.

\vspace{0.5cm}

\paragraph{The copy SB protocol} 
There are two identical sealed labs, indistinguisable from the outside and inside. 
Sleeping Beauty resides in one of those labs.
A fair coin is tossed, with the outcome not revealed to Sleeping Beauty. 
If the outcome is Tails, a copy of Sleeping Beauty is made and put in the second lab, so that in that case both labs contain a Sleeping Beauty. Again Sleeping Beauty knows the protocol.
Each version of Sleeping Beauty is asked whether she thinks the outcome of the coin toss was Heads or Tails, and what her credence for Heads is.
We will consider again the perspective of an outsider Prince Probability. 
All Prince Probability knows is the protocol: he does not know the outcome of the coin toss, whether both labs contain a Sleeping Beauty or only one, which one contains a Sleeping Beauty, \ldots 

\vspace{0.5cm}

\paragraph{Proposed resolution of the copy-SB problem} 
We will argue that if SB is woken up and being asked which outcome the coin toss had, she would assign 1/3 probability for Heads.
An outsider, Prince Probability, who does not know the outcome of the coin toss, gives credence 1/2 to Heads prior to any action he undertakes. 
But when he meets a SB, he changes his credence to 1/3 for Tails. 
Therefore, when Prince Probability and Beauty meet, they agree on their probability assignments for Heads, namely 1/3. Before they have met, however, Prince Probability and Beauty can each correctly assign different credences to Heads, namely 1/2 by Prince Probability and 1/3 by Beauty.

\vspace{0.5cm}

\paragraph{Prince Probability's point of view}
The probability that opening a lab Prince Probability finds Beauty is \begin{equation}
\begin{split}
    p(\text{PP finds SB}) &= p(H) p(\text{PP finds SB}|T)+ p(H) p(\text{PP finds SB}|H) \\ &= \frac{1}{2} \cdot 1 + \frac{1}{2} \frac{1}{2} = \frac{3}{4}    
\end{split}
\end{equation} where `PP finds SB' means the event of Prince Probability finding Sleeping Beauty upon opening the door of a lab. 
The probability of Prince Probability finding a Beauty upon opening the door of a lab under outcome Heads (which we refer to as finding a Heads SB) is 
\begin{equation}
\begin{split}
    p(\text{PP finds Heads SB}) &= p(T) p(\text{PP finds Heads SB}|T)+ p(H) p(\text{PP finds Heads SB}|H) \\ &= \frac{1}{2} \cdot 0 + \frac{1}{2} \cdot \frac{1}{2}   = \frac{1}{4}.    
\end{split}
\end{equation} Therefore, if Prince Probability finds a Sleeping Beauty, the probability that she is a Heads SB (meaning the outcome of the coin toss was Heads) is equal to \begin{equation}
    p(\text{PP finds Heads SB}|\text{PP finds SB}) = \frac{1/4}{3/4} = \frac{1}{3}.
\end{equation}

\vspace{0.5cm}

\paragraph{Interpretation of 1/2 vs 1/3: Prince Probability's point of view.} 
What does it mean that if Prince Probability finds a Beauty, the probability that she is a Heads Beauty is 1/3? Well, before Prince Probability enters a lab, he argues that the probability for Heads is 1/2. 
But, when he enters a lab, then he changes his credence in Heads, depending on whether or not there is a Beauty. 
Namely, if he finds a Beauty, his credence in Heads becomes 1/3. 
This is not surprising: Prince Probability has gotten new information, namely that he sees a Beauty. 
For outcome Tails he is more likely to find a Beauty. Thus, conditioned on Prince Probability finding a Beauty, he gives credence 1/3 to Heads. 
Still, Prince Probability assigns probability 1/2 to Heads before he undertakes any action.

\vspace{0.5cm}

\paragraph{What if Prince Probability finds no Sleeping Beauty?}
If there are only two labs and if Prince Probability opens a lab and finds no Beauty, he will know the outcome was Heads for sure. 
If there were additional empty labs (so under Tails only still two labs contain a Beauty), then Prince Probability will have credence higher than 1/2 but less than 1 that the outcome was Heads, as he is more likely to find no Beauty if the outcome were Heads. His credence here thus depends on the number of labs there are, similar to the discussion of Prince Probability finding Sleeping Beauty asleep in \Cref{sec:SB_Thirder}.

\vspace{0.5cm}

\paragraph{Sleeping Beauty's point of view} 
From the above discussion, we found that when Prince Probability and a Beauty meet, Prince Probability assigns probability 1/3 to Heads. 
But from Beauty's point of view, Prince Probability visiting her does not give her useful information: Prince Probability simply chooses to enter her lab or the other one, with equal probability. 
This seems to suggest that Beauty indeed always assigns probability 1/3.

Consider the following story. God, when he created the universe, tossed a coin Heads or Tails. If the outcome is Tails, there are very favorable conditions, so that more intelligent life is way more likely to exist (i.e. $N$ intelligent organisms are created with $N$ large), whereas for Heads, there are very harsh conditions, making intelligent life harder to exist and reproduce (i.e. only 1 copy exists). 
Now, given that we do exist (and are intelligent, to some extent), which outcome, Heads or Tails, do we give more credence to? 
Intuitively, we would indeed give more credence to Tails. 
This is also our way of doing science: if an explanation makes our observed outcome way more likely, we are more keen to adapt that theory than in the case where an explanation requires extreme fine-tuning (extreme luck, low probability) for our observed outcome.

Similarly, it seems that Sleeping Beauty should consider the event of herself being there as an observation, and incorporate that she is more likely to exist in the case where the outcome is Heads, thus giving more credence to Heads.

Namely, Beauty names the two labs Lab 1 and Lab 2, and without loss of generalisation assumes she is in Lab 2. So the fact that she exists is the event `There is a SB in Lab 2' or shortly `SB in Lab 2'. 
Therefore, her credence in Heads is conditional on the event `SB in Lab 2', an almost trivial observation for her it seems, but one that matters.
Namely, if the outcome of the coin toss is Tails, the probability of `SB in Lab 2' is given by \begin{equation}
    p(\text{SB in Lab 2}|T) = 1, 
\end{equation} whereas for Heads there is a Beauty either in Lab 1 or 2, but not both, \begin{equation}
    p(\text{SB in Lab 2}|H) = \frac{1}{2}.
\end{equation} We also have that \begin{equation}
    p(\text{SB in Lab 2}) = p(T) p(\text{SB in Lab 2}|T) + p(H) p(\text{SB in Lab 2}|H) = \frac{1}{2} \cdot 1 + \frac{1}{2} \cdot \frac{1}{2} = \frac{3}{4}.
\end{equation}
Therefore, the credence for Beauty for Heads is \begin{equation}
    p(T|\text{SB in Lab 2}) = \frac{p(T) p(\text{SB in Lab 2})}{p(\text{SB in Lab 2})} = \frac{1/2 \cdot 1}{3/4} = 2/3.
\end{equation} The credence Sleeping Beauty gives to Heads is thus 1/3, where she uses the fact that she exists (i.e. Lab 2 contains her). 

\vspace{0.5cm}

\paragraph{Summary: Prince Probability's vs Sleeping Beauty's point of view} 
Before opening any lab, Prince Probability assigns credence 1/2 to Heads. However, upon meeting an Sleeping Beauty (opening a lab and seeing a Sleeping Beauty), he changes his assignment to 1/3 for Heads, as he is more likely to meet an Sleeping Beauty in case of Tails. 
Sleeping Beauty, when asked for her credence in Heads, always answers 1/3, and never 1/2. 
Namely, she uses the fact that she exists (i.e. that Lab 2 contains a Sleeping Beauty), which is more likely under Tails.
Therefore, when Prince Probability and Sleeping Beauty meet, their credences for Heads (and Tails) agree. 
Before that, their credences differ: Prince Probability agrees that Sleeping Beauty should assign credence 1/3 to Heads when woken up, but he still correctly assigns credence 1/2 to Heads before opening a lab door.

\subsection{Sleeping Beauty and an unfair coin toss}  \label{sec:unfair_coin}
If the coin toss is unfair, the same calculation as above in \Cref{sec:copy_SB} goes through. 
Suppose that for the coin toss we have $p(H) = 1/4$ and $p(T) = 3/4$, and as before for Tails Sleeping Beauty is copied (so that there are then 2 copies). 
Then, for the case with two labs again\footnote{The same argument goes through if there are more empty labs, and for Tails two of them are filled, whereas for Heads only one is filled.}, we have the following. The probability that opening a lab Prince Probability finds Sleeping Beauty is \begin{equation}
    p(\text{PP finds SB}) = p(H) p(\text{PP finds SB}|H)+ p(T) p(\text{PP finds SB}|T) = \frac{1}{4} \cdot \frac{1}{2} + \frac{3}{4} \cdot 1 = \frac{7}{8}.
\end{equation} The probability of Prince Probability finding a Heads SB upon opening the door of a lab is 
\begin{equation}
    p(\text{PP finds Heads SB}) = p(H) p(\text{PP finds Heads SB}|H)= \frac{1}{4} \cdot \frac{1}{2} = \frac{1}{8}.
\end{equation} Therefore, if Prince Probability finds a SB, the probability that she is a Heads SB (meaning the outcome of the coin toss was Heads) is equal to \begin{equation}
    p(\text{PP finds Heads SB}|\text{PP finds SB}) = \frac{1/8}{7/8} = \frac{1}{7}.
\end{equation} 
Translated to the original Sleeping Beauty problem of \Cref{sec:SB_Thirder}, in the case of un unfair coin toss with $p(H)=1/4$, when Prince Probability finds a Sleeping Beauty awake, he gives credence 1/7 to Heads (and so does Sleeping Beauty too).


\vspace{0.5cm}

\section{Relation to literature on Sleeping Beauty} \label{sec:literature}
Elga's argument using the Principle of Indifference for self-locating belief leads to the same conclusion as ours in \Cref{sec:SB_Thirder}, namely that Sleeping Beauty should assign credence 1/3 to Heads when woken up \cite{elga2020self}. Furthermore, also in all variations considered of the protocol, such as an unfair coin toss and many awakenings, Elga's resolution method gives the same answer as ours, and accords with the updating of the outsider Prince Probability's point of view when meeting a Sleeping Beauty who is awake. 
We can see this fact as an argument for the correctness of Elga's reasoning.

After repeating Elga's argument and showing how it agrees with our arguments, we briefly discuss what arguments on the Sleeping Beauty problem have been given in literature, and counter some of the prominent Halfer arguments. 
We will not be able to mention all relevant results that have been produced about the Sleeping Beauty problem, for which we apologise, but we aim to provide a short overview.

\subsection{Elga's argument for 1/3 and the Principle of Indifference} 
Elga \cite{elga2020self} uses the Principle of Indifference for self-locating belief (see also \cite{elga2004defeating}) to argue that the correct answer for Sleeping Beauty to give when woken up is 1/3. We briefly summarise Elga's argument here, following his writings in \cite{elga2020self} very closely.

\vspace{0.5cm}

If Sleeping Beauty is woken up, she knows she is in one of three events: 
\begin{itemize}
    \item[$H_1$:] Heads and Monday,
    \item[$T_1$:] Tails and Monday,
    \item[$T_2$:] Tails and Tuesday.
\end{itemize}
Namely, upon awakening, Sleeping Beauty clearly assigns zero probability to $H_2$, i.e. the outcome being Heads and it being Tuesday.
So $T = (T_1 \text{ or } T_2)$, $H=(H_1 \text{ or } H_2)$ and $P(H) = P(H_1)$. 

Elga continues as follows: 
\begin{quote}
    `Notice that the difference between your being in $T_1$ and your being in $T_2$ is not a difference in which possible world is actual, but rather a difference in your temporal location within the world. 
    (In a more technical treatment we might adopt a framework similar to the one suggested in Lewis 1983, according to which the elementary alternatives over which your credence is divided are not possible worlds, but rather centered possible worlds: possible worlds each of which is equipped with a designated individual and time. 
    In such a framework, $H_1$, $T_1$, and $T_2$ would be represented by appropriate sets of centered worlds.' (quoted literally from \cite{elga2020self})
\end{quote}

More precisely, in \cite{elga2004defeating} Elga states the Principle of Indifference for self-locating belief as follows (copied from \cite{elga2004defeating}): 
\begin{quote}
    \textsc{Principle of Indifference:} Similar centered worlds deserve equal credence.
\end{quote}

Elga's resolution then proceeds by arguing that $P(T_1) = P(T_2)$ and $P(T_1)=P(H_1)$.

If, when woken up, Sleeping Beauty would be told that the outcome was Tails, she would know she is either in the event $T_1$ or $T_2$. Subjectively, Sleeping Beauty cannot distinguish between $T_1$ and $T_2$, and all the same propositions are True for Sleeping Beauty whether Sleeping Beauty is in $T_1$ or $T_2$. Namely, given Tails, $T_1$ and $T_2$ take place in similar centered worlds, so by the Principle of Indifference for self-locating belief we have $P(T_1|T_1 \text{ or } T_2)=P(T_2|T_1 \text{ or } T_2)$. 

Alternatively to the Sleeping Beauty problem where the coin is flipped \emph{first} on Sunday, we could also imagine the protocol where Sleeping Beauty is woken up on Monday and \emph{then} a coin is flipped determining whether Sleeping Beauty is also woken up on Tuesday (Tails) or not (Heads).
Elga argues that in both cases, if Sleeping Beauty were told that it is Monday, she must assign the same credence to Heads, namely 1/2, i.e. $P(H_1|H_1 \text{ or } T_1) = 1/2$, leading to $P(H_1) = P(T_1)$. 

Thus, using $P(T_1)=P(T_2)$ and $P(H_1)=P(T_1)$ together with the normalisation $P(T_1) + P(T_2) + P(H_1)=1$, we obtain that indeed $P(H)=P(H_1)=1/3$. Elga thus argues that when woken up Sleeping Beauty should give credence 1/3 to Heads, in agreement with our resolution.

\vspace{0.5cm}

Importantly, the Principle of Indifference does \emph{not} allow you to argue that, as $T_1,T_2,H_1$ are all subjectively indistinguishable, they should all get equal credence 1/3. For example, in the case of an unfair coin toss, Elga's argument will correctly divert from credence 1/3 for Heads. For example, if the unfair coin toss gives Heads with probability $h=1/4$, and Tails with probability $1-h=3/4$, then Elga's argument goes as follows. Still, $T_1$ and $T_2$ take place in similar centered worlds, so $P(T_1) = P(T_2)$. Suppose it is Monday, so we condition on the event $H_1 \text{ or } T_1$ (SB waking up on Monday). Then, by the same argument as above, we could also take the coin to be tossed after Sleeping Beauty waking up on Monday and being put back to sleep, for which then Heads has probability $h=1/4$. 
So $P(H_1|H_1 \text{ or } T_1) = h$, giving $P(H_1) = h P(H_1) + h P(T_1)$, and thus $P(H_1) = h/(1-h) P(T_1)$. Combining this with $P(T_1)=P(T_2)$ and the normalisation $P(H_1)+P(T_1)+P(T_2)=1$, we obtain,  for example for $h=1/4$, that $P(H_1)=1/7$. So, in the case of this unfair coin toss, Sleeping Beauty replies giving credence 1/7 to Heads upon being woken up. This accords with our resolution of \Cref{sec:SB_Thirder} and our calculation of the unfair coin toss in \Cref{sec:unfair_coin}.

It is easy to verify that also in the general case of $N$ awakenings and an unfair coin toss with $P(H)=h$, Elga's resolution accords with ours, namely \begin{equation}
    P(H|\text{SB awake}) = \frac{1}{N\cdot \frac{1-h}{h} + 1}.
\end{equation}

\vspace{0.5cm}

As noted by Elga himself using an entertaining story involving Dr Evil \cite{elga2004defeating}, a consequence of the Principle of Indifference that seems somewhat unsettling is the fact that it favors situations where there are many duplicates of ourselves, see also \cite{titelbaum2013ten}. 
Similar difficulties arise when applying the Principle of Indifference to the many-worlds interpretation of quantum mechanics \cite{titelbaum2013ten} and Boltzmann Brain universes\footnote{Another uneasiness about Boltzmann Brains is the following. Boltzmann Brains arise in our models of physics, cosmology and our universe, based on our memories of observed data, but Boltzmann Brains have untrustworthy memories, i.e. their memories of `observed' data are not in line with the universe they live in \cite{carroll2020boltzmann}. Thus if we were Boltzmann Brains, we should not trust our memories of `observed' data, undermining the use of our physical theories and thus the argument for the existence of Boltzmann Brains.} in cosmology \cite{jones2024thinking,titelbaum2013ten}. We will not go deeper into these interesting questions here.

\subsection{The Thirders vs. the Halfers}
Lewis argues in \cite{lewis2001sleeping} against Elga's resolution of 1/3, and favors 1/2, based on the apparent (but false) idea that no new relevant information is gained by Sleeping Beauty waking up. 
In \Cref{sec:SB_Thirder} we have countered this thought by operationalising the SB scenario somewhat more. In that case it is clear that by knowing that she is awake, Beauty does obtain relevant new information, i.e. Beauty being woken up is a nontrivial fact during the two days of the experiment. 
The fact that this piece of information is nontrivial, namely `I (SB) am awake \emph{now}', has also been noted by among others Weintraub, \cite{weintraub2004sleeping} Horgan \cite{horgan2004sleeping}, and Hitchcock \cite{hitchcock2004beauty}, as an argument against 1/2. 
Monton argued for a change in information too, namely that when Beauty wakes up compared to before she was put to sleep on Sunday Beauty has a loss of her temporal location \cite{monton2002sleeping}. 
Moreover, using a variation of the protocol, Dorr \cite{dorr2002sleeping} and Arntzenius \cite{arntzenius2003some} argued how indeed at Beauty's awakening there is nontrivial new information for her.
Our more operational setting and introduction of the outsider's view Prince Probability who may or may not meet Beauty thus provide strong support for these arguments. 

Furthermore, as a consequence of his Halfer position, Lewis \cite{lewis2001sleeping} finds that $P(H|H_1 \text{ or } T_1)=2/3$, i.e. when Beauty wakes up on Monday and learns that it is Monday she should have credence 2/3 in Heads. This is a bit odd, and so-called Double Halfers reject this consequence and instead keep the fact that at all times Beauty should have credence 1/2 in Heads, and should not update her belief. Again this is a bit odd, and Titelbaum \cite{titelbaum2012embarrassment} provided an argument against the Double-Halfer position. 

\vspace{0.5cm}

White provided a generalised Sleeping Beauty problem in \cite{white2006generalized} that, from his thought, seems to argue for 1/2. Namely, a random waking device is used each time Beauty should be woken up, so that only with probability $c \in (0,1]$ she is actually woken up each time the device is used. 
However, for his argument in favour of 1/2 using this generalised Sleeping Beauty problem, White uses a probability conditioned on the event `SB is being woken up at least once', which is different from the self-locating event where Sleeping Beauty notices she has been woken up, also noted in \cite{van2013sleeping}.
Namely, White calculates the probability
\begin{equation}
\begin{split}
    P(H|W) &= \frac{P(H) P(W|H)}{P(H) P(W|H) + P(T) P(W|T} \\ &= \frac{1/2 \cdot c}{1/2 \cdot c + 1/2 \cdot (2c(1-c)+c^2)} = \frac{1}{3-c},    
\end{split}
\end{equation} where $W = \text{(SB awake at least once)}$. White argues that for $c \rightarrow 1$, we obtain $P(H|W)=1/2$, but the conditioned event $W=\text{ `SB awake at least once'}$ is different from `SB awake \emph{now}'.
This is very clear in our discussion of the Sleeping Beauty problem in \Cref{sec:SB_Thirder}. Namely, before Prince Probability enters any lab (and thus assigns credence 1/2 to Heads), he obviously knows that Beauty is woken up at least once. It is only when he actually enters a lab and meets Beauty while she is awake that he updates his credence in Heads to 1/3. 
The fact that for $c \rightarrow 1$, we obtain $P(H|\text{SB awake at least once}) \rightarrow 1/2$ makes sense, as the event `SB awake at least once' in the original protocol for $c=1$ is trivial (has probability 1), and the outsider Prince Probability's credence in Heads, namely 1/2, is recovered.

\vspace{0.5cm}

An argument different from Elga's original argument was given by Dorr and Arntzenius \cite{dorr2002sleeping,arntzenius2003some} using the following variation of the Sleeping Beauty problem. For outcome Heads, Beauty is now awakened as well on Tuesday, but after a short time (while still being awake), she learns that it is Tuesday and the outcome of the coin toss was Heads. Upon waking up, Beauty first assigns equal credence to the four events $T_1,T_2,H_1,H_2$ (Tails Monday, Tails Tuesday, Heads Monday, Heads Tuesday). But if it is not $H_2$, then after a short time, as she does not learn that it is Tuesday and the outcome was Heads, she can eliminate $H_2$. Now, as she learns no information that should influence her credences among the remaining possibilities, still her credences should be divided equally but only among the three events $T_1,T_2,H_1$. Therefore, Beauty should assign credence 1/3 to Heads when being woken up.

\vspace{0.5cm}

Statistical arguments such as Groisman's in \cite{groisman2008end} show that 1/3 and 1/2 answer different questions, and that 1/3 is the correct answer for Beauty `upon the set-up of wakening'. Namely, suppose each time Beauty wakes up, the experimenter who wakes her up and knows the outcome of the coin toss puts a note saying `Heads' or `Tails' in a bowl. After many runs of the protocol, approximately 1/3 of the notes in the bowl will say `Heads'. Furthermore, an investigation of Sleeping Beauty wanting to minimize inaccuracy was investigated by Kierland and Monton in \cite{kierland2005minimizing}, finding that 1/2 and 1/3 may arise from minimizing expected average
inaccuracy or expected total inaccuracy.

The Bayesian rule of updating (conditionalisation) for self-locating beliefs has also been questioned in favour of 1/2 answers for the SB problem in \cite{meacham2008sleeping}. However, having given arguments in favour of 1/3, especially our argument with the/ outsider Prince Probability who when meeting an awake SB changes his credence to 1/3 for Heads questions such proposed change for the Bayesian update rule for self-locating beliefs as in \cite{meacham2008sleeping}.

Finally, there have been Dutch Book arguments \cite{hitchcock2004beauty,bradley2006betting,draper2008diachronic} that make use of rational agents never accepting a bet leading to a sure loss to discuss the Sleeping Beauty problem. 
We will not go into detail here, but refer the reader to \cite{van2013sleeping,titelbaum2013ten}, where an overview of many arguments and discussion accounts on the SB problem can be found.

\vspace{0.5cm}

\section{The Sleeping Beauty problem and the Frauchiger--Renner paradox}  \label{sec:SB_and_FR}
The Frauchiger--Renner (FR) paradox \cite{frauchiger2018quantum} suggests an inconsistency in the naive use of quantum theory to describe the use of itself. In the protocol, agents model each other quantumly and reason about each other's knowledge, arriving at a contradiction. 
The Frauchiger--Renner paradox makes use of a Consistency assumption (C), which can be paraphrased as (copied from \cite{walleghem2024strong}): 
\begin{quote}
    (C): If an agent $B$ knows that statement $s_1$ is true from reasoning with a theory that an agent $A$ accepts and $A$ knows about this reasoning of $B$ about $s_1$, then $s_1$ is true for agent $A$ as well. 
\end{quote}
The extension of the SB problem we provided\footnote{Ref.~\cite{jones2024thinking} provides a similar extension involving such an `outsider'.}, with 1/3 and 1/2 being the correct credences for Sleeping Beauty and an outsider Prince Probability (before meeting SB), respectively, gives an argument against Frauchiger--Renner's assumption (C), as it violates the probabilistic extension thereof \cite{jones2024thinking}.\footnote{In fact, one may perhaps even give an argument for the SB problem refuting (C) (and not only its probabilistic extension). Namely, consider the extension of the SB problem involving the outsider Prince Probability who is asked as well for his belief in Heads. Now we let the experiment take $N$ days, with Tails meaning Beauty is woken up $N$ consecutive days compared to one single day for Tails. Using the same calculus as above in \Cref{sec:SB_Thirder}, Beauty has credence $1/N$ in Heads, whereas Prince Probability has credence 1/2 in Heads. Taking $N \rightarrow \infty$ Beauty now has credence 0 in Heads, whereas Prince Probability has non-zero credence 1/2 in Heads. I thank Rui Soares Barbosa for this note.}
However, this does not resolve the stronger GHZ--FR paradox of \cite{walleghem2024strong} that is similar to the FR one, but makes use of a weaker Consistency assumption where only classically communicating classical agents (who may physically meet and greet each other if desired) need to agree on statements. 
This assumption is not violated in the SB problem or its extension, as indeed Sleeping Beauty and the outsider Prince Probability agree when they meet each other. 
This suggests that these FR-like paradoxes in quantum theory need a more radical resolution rather than refuting the Consistency assumption (C).

\vspace{0.5cm}

\section{Discussion and conclusion} \label{sec:conclusion}
In this short note we have given arguments employing simple Bayesian probability that suggest 1/3 to be the correct answer for Sleeping Beauty to reply when she is woken up and asked for her credence in Heads. 
Sleeping Beauty does have extra useful information, which is why 1/2 is not her correct answer. 
An outsider, Prince Probability, assigns credence 1/2 to Heads, but when Prince Probability meets Sleeping Beauty being awake, he lowers his credence in Heads to 1/3.
Furthermore, if Prince Probability meets Sleeping Beauty while asleep, he instead increases his credence in Heads.

\vspace{1cm}

\section*{Acknowledgements} 
I am most grateful to Caroline Jones for bringing this problem to my attention, and for many following enlightening discussions that lead to this short note.
I also thank Will Oliver and Rui Soares Barbosa for discussions on the statistical interpretation of the problem. I thank Khadija Sarr for the name Prince Probability and the title.

LW acknowledges support from the United Kingdom Engineering and Physical Sciences Research Council (EPSRC) DTP Studentship (grant number EP/W524657/1). LW also thanks the International Iberian Nanotechnology Laboratory -- INL in Braga, Portugal and the Quantum and Linear-Optical Computation (QLOC) group for the kind hospitality.

\newpage

\printbibliography

@article{jones2024thinking,
  title={Thinking twice inside the box: is {W}igner's friend really quantum?},
  author={Jones, Caroline L and Müller, Markus P},
  year={2024},
journal = {arXiv preprint},
doi = {10.48550/arXiv.2402.08727}
}

@article{groisman2008end,
  title={The end of {S}leeping {B}eauty's nightmare},
  author={Groisman, Berry},
  journal={The British Journal for the Philosophy of Science},
    number = {3},
    pages = {409--416},
    volume = {59},
    year={2008},
    publisher = {Oxford University Press},
    doi = {10.1093/bjps/axn015}
}

@article{elga2020self,
  title={Self-locating belief and the {S}leeping {B}eauty problem},
  author={Elga, Adam},
  journal={Analysis},
  Volume={60},
number = {2},
    doi = {10.1093/analys/60.2.143},
  year={2000},
pages = {143-147}
}

@article{lewis2001sleeping,
  title={Sleeping {B}eauty: reply to {E}lga},
  author={Lewis, David},
  journal={Analysis},
  volume={61},
  number={3},
  pages={171--176},
  year={2001},
  publisher={JSTOR},
    doi = {10.1093/analys/61.3.171} 
}

@article{elga2004defeating,
  title={Defeating {D}r. {E}vil with self-locating belief},
  author={Elga, Adam},
  journal={Philosophy and Phenomenological Research},
  volume={69},
  number={2},
  pages={383--396},
  year={2004},
  publisher={Wiley Online Library},
doi ={10.1111/j.1933-1592.2004.tb00400.x}
}

@article{meacham2008sleeping,
  title={{S}leeping {B}eauty and the dynamics of de se beliefs},
  author={Meacham, Christopher JG},
  journal={Philosophical Studies},
  volume={138},
number = {2},
  pages={245--269},
  year={2008},
  publisher={Springer},
	doi = {10.1007/s11098-006-9036-1}
}

@article{weintraub2004sleeping,
  title={{S}leeping {B}eauty: a simple solution},
  author={Weintraub, Ruth},
  journal={Analysis},
  volume={64},
  number={1},
  pages={8--10},
  year={2004},
  publisher={JSTOR},
    doi = {10.1093/analys/64.1.8}
}

@article{horgan2004sleeping,
  title={{S}leeping {B}eauty awakened: New odds at the dawn of the new day},
  author={Horgan, Terry},
  journal={Analysis},
  volume={64},
  number={1},
  pages={10--21},
  year={2004},
    doi = {10.1093/analys/64.1.10}
}

@article{white2006generalized,
  title={The generalized {S}leeping {B}eauty problem: a challenge for thirders},
  author={White, Roger},
  journal={Analysis},
  volume={66},
  number={2},
  pages={114--119},
  year={2006},
  publisher={JSTOR},
    doi = {10.1093/analys/66.2.114}
}

@mastersthesis{van2013sleeping,
  title={On the {S}leeping {B}eauty problem},
  author={van der Meijden, SP},
  type={{B.S.} thesis},
  year={2013},
school = {Utrecht University},
url = {https://studenttheses.uu.nl/handle/20.500.12932/14598}
}

@article{dorr2002sleeping,
  title={{S}leeping {B}eauty: in defence of {E}lga},
  author={Dorr, Cian},
  journal={Analysis},
  volume={62},
  number={4},
  pages={292--296},
  year={2002},
doi = {10.1093/analys/62.4.292}
}

@book{titelbaum2013quitting,
  title={\href{https://doi.org/10.1017/S0266267114000492}{Quitting certainties: A {B}ayesian framework modeling degrees of belief}},
  author={Titelbaum, Michael G},
  year={2013},
  publisher={Oxford University Press},
doi = {https://doi.org/10.1017/S0266267114000492}
}

@article{titelbaum2012embarrassment,
  title={An embarrassment for double-halfers},
  author={Titelbaum, Michael G},
  journal={Thought: A Journal of Philosophy},
  volume={1},
  number={2},
  pages={146--151},
  year={2012},
  publisher={Wiley Online Library},
	doi = {10.1002/tht3.21}
}

@article{arntzenius2003some,
  title={Some problems for conditionalization and reflection},
  author={Arntzenius, Frank},
  journal={The Journal of Philosophy},
  volume={100},
  number={7},
  pages={356--370},
  year={2003},
  publisher={JSTOR},
url ={http://www.jstor.org/stable/3655783}
}

@article{monton2002sleeping,
  title={{S}leeping {B}eauty and the forgetful {B}ayesian},
  author={Monton, Bradley},
  journal={Analysis},
  volume={62},
  number={1},
  pages={47--53},
  year={2002},
  doi = {10.1093/analys/62.1.47}
}

@article{hitchcock2004beauty,
  title={Beauty and the bets},
  author={Hitchcock, Christopher},
  journal={Synthese},
  volume={139},
  pages={405--420},
  year={2004},
  publisher={Springer},
doi = {https://doi.org/10.1023/B:SYNT.0000024889.29125.c0}
}

@article{piccione1997interpretation,
  title={On the interpretation of decision problems with imperfect recall},
  author={Piccione, Michele and Rubinstein, Ariel},
  journal={Games and Economic Behavior},
  volume={20},
  number={1},
  pages={3--24},
  year={1997},
  publisher={Elsevier},
issn = {0899-8256},
doi = {https://doi.org/10.1006/game.1997.0536}
}

@article{aumann1997forgetful,
  title={The forgetful passenger},
  author={Aumann, Robert J and Hart, Sergiu and Perry, Motty},
  journal={Games and Economic Behavior},
  volume={20},
  number={1},
pages = {117-120},
year = {1997},
issn = {0899-8256},
doi = {10.1006/game.1997.0578}
}

@article{bradley2006betting,
    title={When betting odds and credences come apart: More worries for {D}utch book arguments},
    author={Bradley, Darren and Leitgeb, Hannes},
    journal={Analysis},
    volume = {66},
    number = {2},
    pages = {119-127},
    year = {2006},
    month = {04},
    issn = {0003-2638},
    doi = {10.1093/analys/66.2.119}
}

@article{draper2008diachronic,
  title={Diachronic {D}utch books and {S}leeping {B}eauty},
  author={Draper, Kai and Pust, Joel},
  journal={Synthese},
  volume={164},
  pages={281--287},
  year={2008},
  publisher={Springer},
doi = {https://doi.org/10.1007/s11229-007-9226-1}
}

@article{titelbaum2013ten,
  title={Ten reasons to care about the {S}leeping {B}eauty problem},
  author={Titelbaum, Michael G},
  journal={Philosophy Compass},
  volume={8},
  number={11},
  pages={1003--1017},
  year={2013},
 publisher = {Wiley-Blackwell},
doi = {10.1111/phc3.12080}
}

@incollection{carroll2020boltzmann,
  title={Why {B}oltzmann brains are bad},
  author={Carroll, Sean M},
  booktitle={Current controversies in philosophy of science},
  pages={7--20},
  year={2020},
  publisher={Routledge},
 url = {https://doi.org/10.48550/arXiv.1702.00850}
}

@article{kierland2005minimizing,
	author = {Brian Kierland and Bradley Monton},
	doi = {10.1111/j.1933-1592.2005.tb00533.x},
	journal = {Philosophy and Phenomenological Research},
	number = {2},
	pages = {384--395},
	publisher = {Wiley-Blackwell},
	title = {Minimizing Inaccuracy for Self-Locating Beliefs},
	volume = {70},
	year = {2005}
}

@misc{neal2006puzzles,
  title={\href{https://doi.org/10.48550/arXiv.math/0608592}{Puzzles of anthropic reasoning resolved using full non-indexical conditioning}},
  author={Neal, Radford M},
  year={2006},
 eprint = {math/0608592},
archivePrefix = {arXiv},
primaryclass = {math.ST}
}

@book{stalnaker2010our,
  title={\href{https://doi.org/10.1093/acprof:oso/9780199545995.001.0001}{Our knowledge of the internal world}},
  author={Stalnaker, Robert},
  year={2008},
month = {07},
  publisher={Oxford University Press},
    isbn = {9780199545995},
    doi = {10.1093/acprof:oso/9780199545995.001.0001},
    url = {https://doi.org/10.1093/acprof:oso/9780199545995.001.0001}
}

@article{walleghem2024strong,
    title = {A refined {F}rauchiger--{R}enner paradox based on strong contextuality},
    author = {Walleghem, Laurens and Barbosa, Rui Soares and Pusey, Matt and Weigert, Stefan},
year = {2024},
journal = {to appear}
}

@article{frauchiger2018quantum,
  title={Quantum theory cannot consistently describe the use of itself},
  author={Frauchiger, Daniela and Renner, Renato},
  journal={Nature communications},
  volume={9},
  number={1},
  pages={3711},
  year={2018},
  publisher={Nature Publishing Group UK London},
  doi={https://doi.org/10.1038/s41467-018-05739-8}
}
\end{document}